\documentclass{elsart}
\usepackage{graphicx}
\usepackage{psfrag}
\usepackage{times}
\usepackage{amsbsy, amsmath}

\newcommand{\pd}{\partial}

\newcommand{\bu}{{\mathbf u}}

\newcommand{\bB}{{\mathbf B}}

\newcommand{\bx}{{\mathbf x}}
\newcommand{\bxp}{{\mathbf x'}}
\newcommand{\xp}{x^\prime}

\newcommand{\lap}{\Delta}

\newcommand{\grad}{{\nabla}}

\newcommand{\curl}{{\nabla} \times}
\newcommand{\divv}{{\nabla} \cdot}

\newcommand{\ey}{\hat{\mathbf e}_y}

\newcommand{\mT}{\mathcal{T}}

\newcommand{\bdy}{{\partial\Omega}}

\newcommand{\pdn}{\partial_n}

\newcommand{\ex}{\hat{\mathbf e}_x}

\input epsf
\usepackage{amssymb,amsmath}
\begin{document}
\begin{frontmatter}
\title{Spectral method for matching exterior and interior elliptic problems}
\author{Piotr Boronski}


\address{Center for Turbulence Research\\
Stanford University, Bldg. 500, Stanford, CA 94305-3035 \\
e-mail: \underline{boronski@gmail.com}}


\begin{abstract}
A spectral method is described for solving coupled elliptic problems on an
interior and an exterior domain.  The method is formulated and tested on the
two-dimensional interior Poisson and exterior Laplace problems, whose
solutions and their normal derivatives are required to be continuous across
the interface.  A complete basis of homogeneous solutions for the interior and
exterior regions, corresponding to all possible Dirichlet boundary values at
the interface, are calculated in a preprocessing step.  This basis is used to
construct the influence matrix which serves to transform the coupled boundary
conditions into conditions on the interior problem.  Chebyshev approximations
are used to represent both the interior solutions and the boundary values.  A
standard Chebyshev spectral method is used to calculate the interior
solutions.  The exterior harmonic solutions are calculated as the convolution
of the free-space Green's function with a surface density; this surface
density is itself the solution to an integral equation which has an analytic
solution when the boundary values are given as a Chebyshev expansion.
Properties of Chebyshev approximations insure that the basis of exterior
harmonic functions represents the external near-boundary solutions uniformly.
The method is tested by calculating the electrostatic potential resulting from
charge distributions in a rectangle.  The resulting influence matrix is
well-conditioned and solutions converge exponentially as the resolution is
increased.  The generalization of this approach to three-dimensional problems
is discussed, in particular the magnetohydrodynamic equations in a finite
cylindrical domain surrounded by a vacuum.
\end{abstract}

\end{frontmatter}

\section{Motivation}
\label{sec:motivation}
The search for a self-sustaining magnetohydrodyamic dynamo has taken on 
great momentum in recent years, as researchers have sought to produce dynamos 
in the laboratory \cite{Steglitz02,Riga01,VKS02,Forest02,Shew02} and
in simulations \cite{Jepps67,Dudley89,Glatzmaier95,Tilgner97,Hollerbach00,Willis02,Matsui04,Iskakov04,Xu04,Guermond06}.
One of the fundamental problems in numerical magnetohydrodynamics is the
formulation of boundary conditions. The governing equations 
describe the velocity and magnetic field in
in a finite container of electrically conducting fluid.
At the container boundaries, the velocity is specified,
but the magnetic field is not.
Instead, the magnetic field is required to satisfy
continuity conditions with the exterior magnetic field
in the domain surrounding the fluid.
The nature of these conditions depends on the properties
of the surrounding medium; a complete discussion can be 
found in \cite{Roberts67}.

Under the quasi-static approximation \cite{Roberts67},
for a given velocity field $\bu$ and magnetic Reynolds number $Rm$,
the equations describing the interior magnetic field are 
\begin{subequations}
\label{eq:Bint}
\begin{eqnarray}
\pd_t \bB &=& \curl(\bu\times\bB)+ \frac{1}{Rm} \lap \bB \label{eq:induction}\\
\divv \bB &=& 0 \label{eq:nomonopoles}
\end{eqnarray}
\end{subequations}
The case of a fluid of finite electric conductivity restricted to a
finite volume and surrounded by vacuum is of special
importance because it models a number of experimental, geophysical,
and astrophysical configurations.
Since there are no electrical currents in a vacuum, $\bB^{vac}$ is
curl-free, and is therefore the gradient of a potential 
if the exterior domain is simply connected.
The exterior magnetic field then obeys:
\begin{subequations}
\begin{eqnarray}
\bB^{\rm vac}&=&\grad\phi^{\rm vac}\\
\lap\phi^{\rm vac} &=& 0  \label{eq:phiharm}\\
\bB^{\rm vac} &\rightarrow& 0\qquad |\bx|\rightarrow\infty
\end{eqnarray}
\label{eq:Bext}
\end{subequations}
but is otherwise not fixed. 
The magnetic field is required to be continuous at the boundary:
\begin{equation}
\bB - \bB^{\rm vac} = 0 \qquad \bx\in\bdy 
\label{eq:Bcont}
\end{equation}
In this case,
continuity of all three components of the magnetic field are sufficient
to uniquely determine both the interior and exterior fields. 
Our ultimate goal is to transform \eqref{eq:Bext}-\eqref{eq:Bcont}
into boundary conditions that can be applied to \eqref{eq:Bint}
without calculating $\bB^{\rm vac}$.
The general principle we will employ is to construct a complete basis of
exterior solutions $\bB^{\rm vac}$ of \eqref{eq:Bext} in a preprocessing step, 
and to calculate $\bB|_\bdy$ for each member of the basis.
The matching conditions \eqref{eq:Bcont} will then yield boundary conditions
for $\bB$.

To explore this approach, we will apply it to the simpler analogous
scalar problem of the Poisson problem in an interior domain.
We will require the solution to match continuously to an exterior
solution satisfying Laplace's equation.
In this case, both Dirichlet and Neumann matching conditions are 
necessary to specify a unique solution.
Formally, we wish to solve the following problem:
\begin{subequations}
\label{eq:two_domains}
\begin{eqnarray}
\lap \Phi=\rho  &&\quad\text{in }\Omega \\
\lap \phi=0 && \quad \text{outside }\Omega
\end{eqnarray}
\end{subequations}
with boundary conditions:
\begin{subequations}
\label{eq:match}
\begin{eqnarray}
\Phi(\bx)-\phi(\bx) &=& 0 \qquad \bx\in\bdy \label{eq:matchD}\\
\pdn\Phi(\bx)-\pdn\phi(\bx) &=& 0 \qquad \bx\in\bdy \label{eq:matchN}\\
\nabla\phi(\bx)&\rightarrow& 0 \qquad |\bx|\rightarrow\infty
\end{eqnarray}
\end{subequations}
where $\Omega$ is a bounded domain with boundary $\bdy$.
A physical interpretation of \eqref{eq:two_domains}-\eqref{eq:match} 
is that of an electrostatic potential $\Phi$ of a field 
generated by charges distributed in space with the density $-\rho$,
where the electrical permeability of the vacuum is taken to be one
by the choice of units.
We wish to calculate the interior solution $\Phi$ 
without explicitly constructing $\phi$.

There exists a vast literature on the numerical solution of the
fundamental physical problems \eqref{eq:two_domains}--\eqref{eq:match} 
and \eqref{eq:Bint}-\eqref{eq:Bcont}. We will briefly
survey a small portion of this literature here, and postpone 
a more detailed comparison between our method and others
to a later section.

The main tool by which exterior domains can be eliminated is
Green's theorem, which replaces elliptic differential equations 
over a domain with integrals over the bounding surfaces.
The use of methods based on boundary integrals has grown explosively 
since the 1970s-1980s to solve engineering problems from fields such as 
acoustics, elasticity, electromagnetism and fluid mechanics
\cite{Hsiao73,Jaswon77,Brebbia78,Banerjee81,Atkinson82,Rokhlin83,Brebbia83,Pozrikidis92,KirkupWeb}.
Hybrid methods, coupling a differential equation formulation in
a domain and a boundary integral formulation at the boundary
via an influence matrix, were also developed at the same time
to solve \eqref{eq:two_domains}--\eqref{eq:match} and similar problems.
The majority of these approaches have been based on finite elements
and are hence applicable to complicated real-world geometries.
The boundary integrals are discretized with techniques derived 
from finite element theory, leading to the term boundary element method,
and the hybrid methods use finite elements to solve the equations 
in the domain.

This situation contrasts with magnetohydrodynamics, which has
been dominated by spectral methods.
Spherical domains are standard, for geophysical and astrophysical reasons.
Spectral methods can then be based on spherical harmonics and the 
poloidal-toroidal decomposition \cite{Dudley89,Glatzmaier95,Tilgner97,Hollerbach00}.
The solution to Laplace's equation on the exterior is immediate and,
moreover, solutions and associated boundary conditions for
each spherical harmonic and toroidal or poloidal component are decoupled.
Boundary conditions at the interface
can then be formulated for each mode without the use of an influence matrix.

The technique which we will describe is based on spectral methods,
but the geometry is assumed to be somewhat more complicated.
Our technique occupies the niche which spectral methods staked out 
in the 1970s-1980s when the use of Chebyshev polynomials became 
commonplace \cite{Orszag71,Canuto88} to represent domains with 
one or more non-periodic directions.
In keeping with this tradition, we expect its main application to be 
to tensor-product domains whose boundaries consist of a small number of 
piecewise-smooth surfaces, such as the finite three-dimensional 
cylinder which is our eventual goal.

We also mention here some other recent approaches
to solving the magnetohydrodynamic equations
\cite{Matsui04,Iskakov04,Iskakov05,Xu04,Guermond06},
with a view to generalizing the geometry and/or increasing parallelization.
In \cite{Iskakov04,Iskakov05}, a finite volume method is
used to discretize the solution in the interior, which is
matched to that in the exterior vacuum via a boundary element method.
\cite{Xu04} describes an integral equation formulation for the
entire domain, and \cite{Guermond06} uses finite elements with
a penalty method to apply boundary conditions.

\section{Influence matrix formulation}
\label{sec:Influence matrix formulation}

We formulate a two-stage method for 
solving \eqref{eq:two_domains}--\eqref{eq:match},
consisting of an initial preprocessing step which depends only on
the geometry, followed by a step whereby solutions for many 
different distributions $\rho$ can be generated at 
little incremental cost.
This is the usual description of the decomposition of $\Phi$ into 
homogeneous and particular solutions, with the additional
proviso that solutions in the exterior domain are to be taken
into account in the preparation of the homogeneous solutions.
We will construct the homogeneous solutions by generating harmonic bases
$\{\Phi_j^h\}$ and $\{\phi_j\}$ of interior and exterior solutions,
corresponding to Dirichlet boundary data $\{f_j\}$ to be specified later.
We decompose \eqref{eq:two_domains}--\eqref{eq:match} 
into the Poisson and Laplace problems:
\begin{equation}
\begin{array}{lll}
\lap\Phi^p=\rho \;\;\text{in } \Omega\qquad & 
\qquad \lap\Phi_j^h=0 \;\;\text{in } \Omega\qquad & 
\qquad \lap\phi_j = 0 \;\;\text{outside }\Omega\\
\Phi^p|_\bdy=0 \qquad & \qquad \Phi_j^h|_\bdy = f_j \qquad & \qquad \phi_j|_\bdy = f_j\\
&&\qquad\nabla\phi_j|_\infty = 0 
\end{array}
\label{eq:construct}\end{equation}
and then construct the linear superpositions:
\begin{equation}
\begin{array}{ccc}
\Phi=\Phi^p+\Phi^h  \qquad & \qquad \Phi^h =\sum_j c_j\Phi_j^h \qquad & \qquad \phi=\sum_j c_j\phi_j
\end{array}
\label{eq:super}\end{equation}
Then
\begin{eqnarray}
\lap\Phi = &\lap\Phi^p+\sum_j c_j\lap\Phi_j^h = \rho &\qquad\text{ in } \Omega\\
\lap\phi = &\sum_j c_j\lap\phi_j = 0 &\qquad\text{ outside }\Omega\\
\left.(\Phi-\phi)\right|_\bdy = &\Phi^p|_\bdy +
\sum_j c_j \left.(\Phi^h_j-\phi_j)\right|_\bdy = 0 & \\
\nabla\phi|_\infty =& \sum_j c_j \:\nabla\phi_j|_\infty = 0 &
\end{eqnarray}
are already satisfied by construction, while
\begin{equation}
\left.\pdn(\Phi-\phi)\right|_\bdy
= \left.\pdn{\Phi^p}\right|_\bdy + 
\sum_j c_j\left.\pdn(\Phi^h_j-\phi_j)\right|_\bdy = 0
\label{eq:matchN2}\end{equation}
constitutes a system of equations to be solved for ${c_j}$,
where the derivative with respect to the normal is taken 
in the direction from the interior to the exterior region
for both $\Phi$ and $\phi$.
$\Phi$ is then set equal to the sum in \eqref{eq:super}.
If the interior harmonic functions are not stored, $\Phi$ can
be obtained by solving:
\begin{eqnarray}
\lap\Phi &=& \rho  \;\;\text{in } \Omega \\
\Phi|_\bdy &=& \sum_j c_j f_j
\end{eqnarray}
Using $\bx_i$ to index points on the boundary, 
\eqref{eq:matchN2} can be discretized as:
\begin{equation}
\sum_j \left[\pdn(\phi_j-\Phi^h_j)(\bx_i) \right]c_j = \left.\pdn{\Phi^p}\right(\bx_i)
\label{eq:matchN3}
\end{equation}
Equation \eqref{eq:matchN3} shows that the goal of the preprocessing
step is the construction and inversion of
the influence or capacitance matrix: 
\begin{equation}
C_{ij} \equiv \left[\pdn\left(\phi_j - \Phi^h_j\right)\left({\bx_i}\right) \right]
\label{eq:def_influ_matrix}
\end{equation}
The functions $\{f_j\}$ are required to constitute a complete set
for values along the discretized boundary $\bdy$. 
Another way to describe the influence matrix is as a discrete
representation of the difference between the 
Dirichlet-to-Neumann mappings in the exterior and in the interior regions.

Equivalently, $\pdn(\phi_j-\Phi^h_j)|_\bdy$ can be represented as
coefficients of a basis set $\{g_i\}$ (which may be identical with 
the set of boundary value functions $\{f_i\}$) along each boundary.
Equation \eqref{eq:matchN2} is then discretized as:
\begin{equation}
\sum_j \langle\pdn(\phi_j-\Phi^h_j),g_i\rangle \: c_j = 
\langle \pdn{\Phi^p},g_i \rangle
\label{eq:matchN4}
\end{equation}
Although we will use $\bx_i$ and the notation in \eqref{eq:matchN3}
in what follows, the method is easily reformulated using \eqref{eq:matchN4}.

\section{Solution of Poisson and Laplace problems}

\subsection{Interior domain}
\label{sec:Interior domain}

We now turn to the solution of \eqref{eq:construct}.
For the interior problems listed in the first two columns,
we assume that we dispose of a solver able to
compute solutions to Poisson's equation in $\Omega$ 
with any specified boundary values. 
In principle, any numerical method can be used.
In our particular case, we use a spectral discretization \cite{Canuto88}
\begin{equation}
\Phi(x,y) = \sum_{k,l=0}^{K,L}\mT_k(x/H)\mT_l(y)
\label{eq:ChebCheb}\end{equation}
for the rectangle $[-H,H]\times[-1,1]$. 
The spectral basis functions are the Chebyshev polynomials 
$\mT_k(x)=\cos(k\arccos(x))$. Taking $H\geq 1$, we set $K\geq L$.
We use a standard method \cite{Canuto88} to solve the Poisson equation with 
Dirichlet boundary conditions, diagonalizing the discretized 
second derivative operator in $y$, and using recursion relations to
treat the second derivative in $x$.

\subsection{Exterior harmonic functions}
\label{sec:Exterior harmonic functions}

Our main focus is on the construction of the exterior harmonic solutions,
specified in the third column of \eqref{eq:construct}. 
In order to avoid truncating or spatially discretizing the exterior domain,
we will construct $\{\phi_j\}$ using the fundamental 
solution of the Laplace equation: the Green's function satisfying
\begin{eqnarray}
\lap_\bxp G(\bx ; \bxp) =& \delta(\bx-\bxp) & \\
\nabla G(\bx ; \bxp) =& 0  & \qquad\text{for } \bxp \rightarrow \infty
\label{eq:fundamental_Green}\end{eqnarray}
For a specified boundary value distribution $f(\bx)$, 
we first calculate an appropriate source distribution $\sigma(\bx)$ on the 
boundary by solving the integral equation:
\begin{equation}
\int_\bdy G(\bx ; \bxp) \sigma (\bxp) = f(\bx) \qquad\text{for }\bx\in\bdy
\label{eq:chargedist}\end{equation}
The exterior harmonic function $\phi(\bx)$ required is then:
\begin{equation}
\phi(\bx) \equiv \int_\bdy G(\bx ; \bxp) \sigma (\bxp)
\label{eq:phiconstruct}\end{equation}
where $\bx$ takes values either on or off $\bdy$.

We now apply \eqref{eq:chargedist}-\eqref{eq:phiconstruct}
to our particular test problem of a rectangle.
We divide the set of boundary distributions into four sets, 
each taking non-zero values on only one side of the rectangle. 
In two dimensions, the fundamental Green's function solving
\eqref{eq:fundamental_Green} is
\begin{equation}
\frac{-1\;}{2\pi}\ln|\bx-\bxp|
\end{equation}
Equation \eqref{eq:chargedist} thus reduces to:
\begin{equation}
	\int_a^b \frac{-1\;}{2\pi}\ln|x-\xp|\sigma(\xp)\,d\xp = f(x)
	\label{eq:Carleman_equation}
\end{equation}
Equation \eqref{eq:Carleman_equation} is known as 
Symm's or Carleman's equation 
and has the following solution \cite{Jorgens70,Polyanin98}:
\begin{eqnarray}
\sigma(x)=\frac{-2}{\pi\sqrt{(x-a)(b-x)}}
\left[\int_a^b\frac{\sqrt{(\xp-a)(b-\xp)}f'(\xp)\,d\xp}{\xp-x}\right. &&\nonumber\\
+\frac{1}{\ln ((b-a)/4)}\left.\int_a^b\frac{f(\xp)\,d\xp}{\sqrt{(\xp-a)(b-\xp)}}\right]&&
\label{eq:Carleman_sol}
\end{eqnarray}
if $b-a\neq 4$. 
(If $b-a=4$ then the second integral in \eqref{eq:Carleman_sol}
can be replaced by an arbitrary constant if 
$\int _a^b f(t)\left[(t-a)(b-t)\right]^{-\frac{1}{2}}\,dt=0$.)

Up to now, we have not specified the Dirichlet boundary values $f$.
The choice of boundary value distributions is restricted 
only by the requirement, stated in section 
\ref{sec:Influence matrix formulation}, 
that the set of distributions form a basis for functions defined 
on the boundary $\bdy$.
Because we use Chebyshev polynomials to represent the interior solutions,
it is convenient to take as boundary values $f_k(x)$ each of the functions
$\mT_k(x/H)$ on the interval $[-H,H]$.
The corresponding solutions $\sigma_k(x;H)$ obtained from evaluating
\eqref{eq:Carleman_sol} are:
\begin{equation}
\sigma_k(x;H)= A_k\frac{\mT_k(x/H)}
{\pi H\sqrt{1-\left(\frac{x}{H}\right)^2}} \quad;\quad 
A_k=\left\{\begin{array}{ll}
2\pi k\quad & k>0\\
-2\pi\left[\ln(H/2)\right]^{-1}\quad & k=0
\end{array}\right.
\label{eq:Carleman_sol_Cheb}\end{equation}

This remarkable property -- the fact that
that weighted Chebyshev polynomials are also obtained as
the source distributions $\sigma_k(x)$ when the boundary values
$f_k(x)$ are Chebyshev polynomials -- is related to the very
reason that Chebyshev polynomials are optimal in approximating
polynomials on the interval.
The function $1/\pi\sqrt{1-x^2}$ in \eqref{eq:Carleman_sol_Cheb} (for $H=1$)
is the weight with respect to which Chebyshev polynomials are orthogonal 
on the interval and is the asymptotic density of the Chebyshev
interpolation points $\cos(\pi j/J)$, the extrema of the Chebyshev
polynomials. See \cite{Mason00,Trefethen00} for further details.
Note also that the orthogonality of the Chebyshev
polynomials with respect to this weight causes the
second integral in \eqref{eq:Carleman_sol} to vanish except for $\mT_0$.

\begin{figure}[htbp]
{\centering
\includegraphics[width=8cm]
{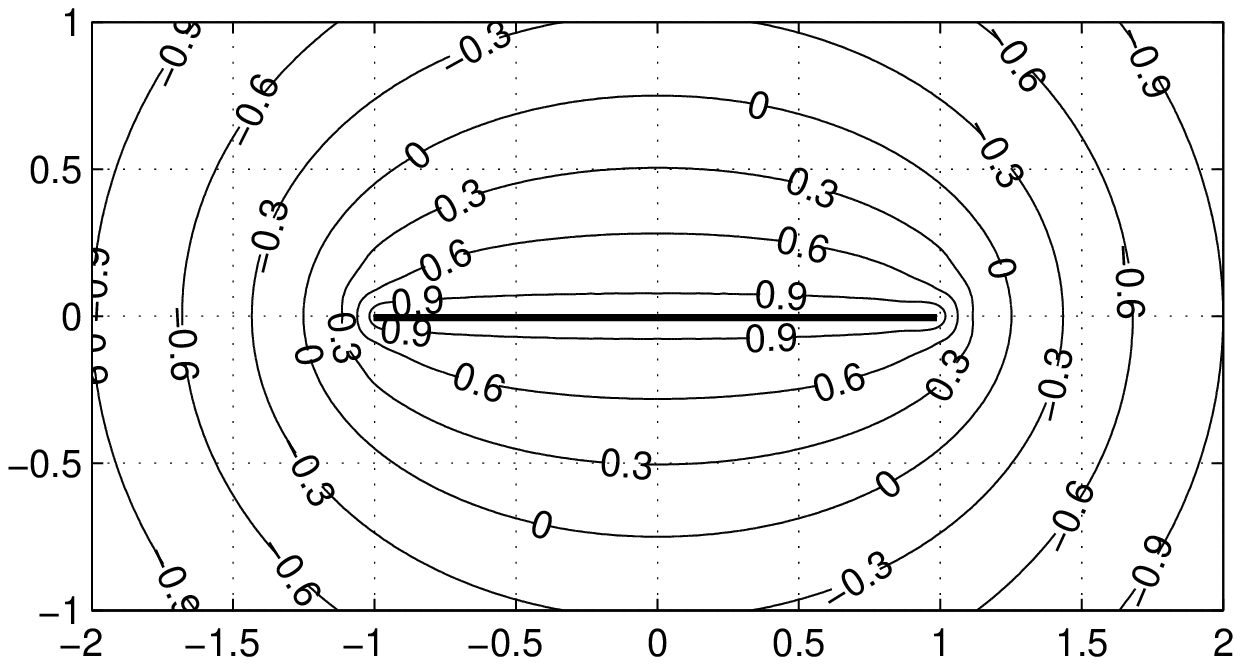}

\includegraphics[width=8cm]
{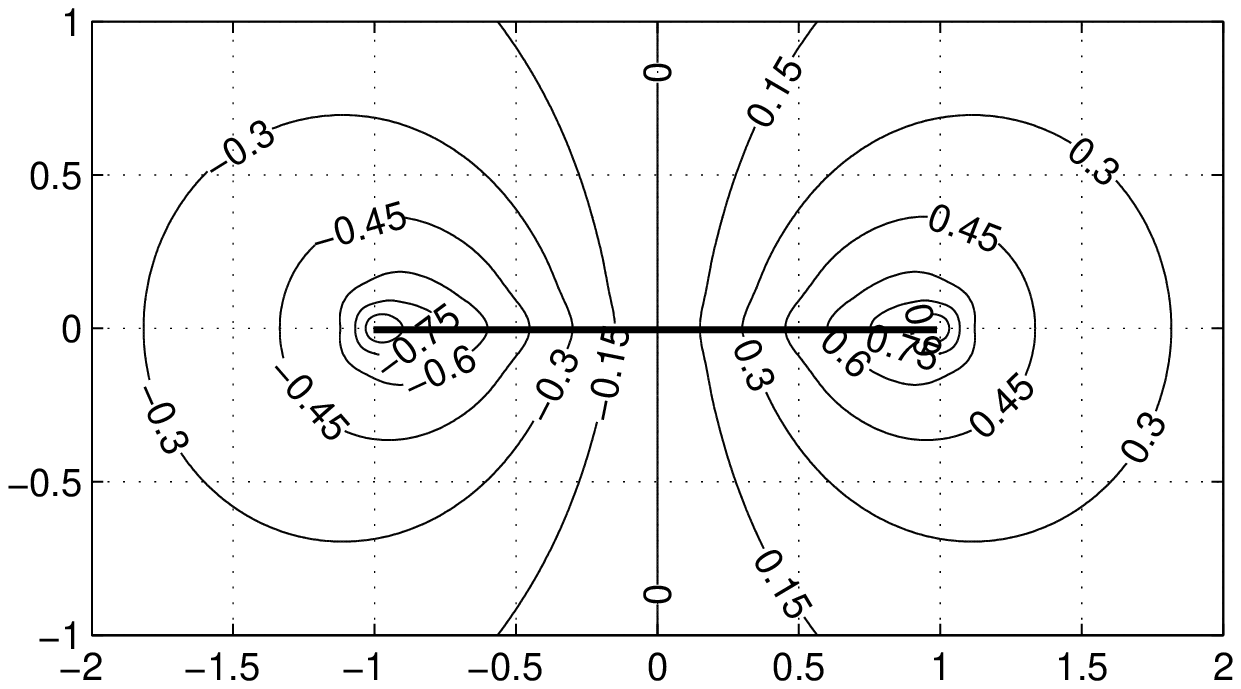}

\includegraphics[width=8cm]
{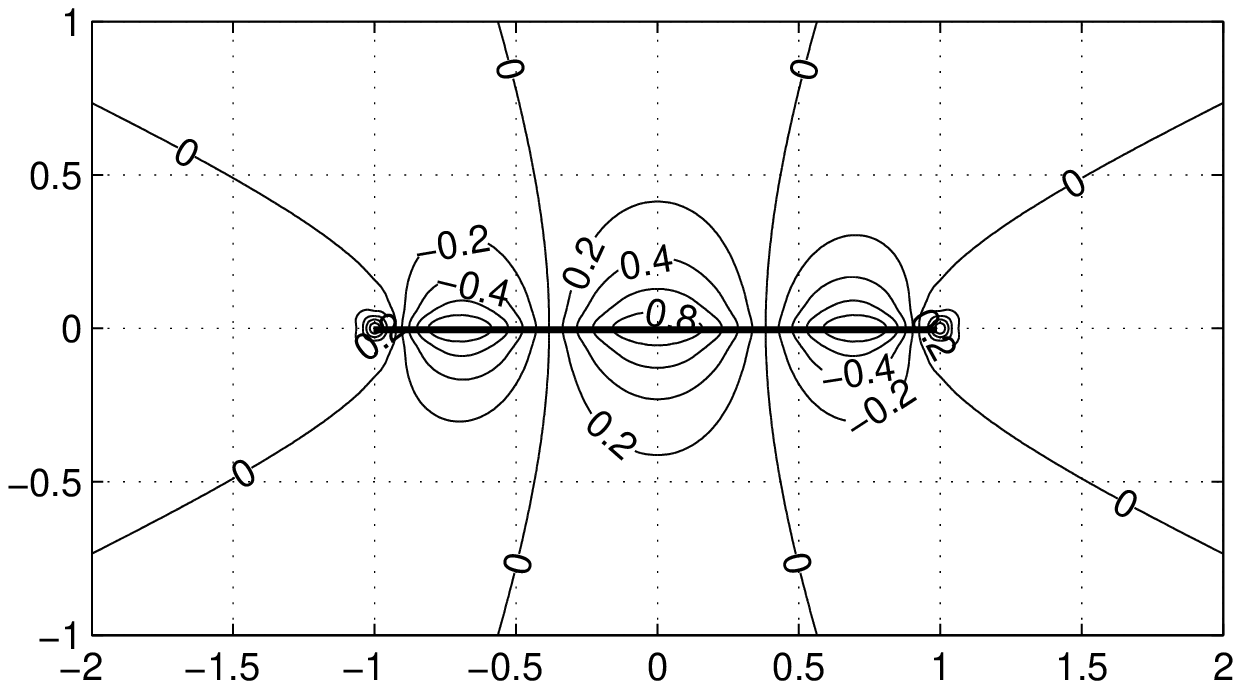}

}
\caption{Potentials $\phi^x_k(\bx)$ generated by line source distributions
$\sigma_k(x;H=2)$. Values of the potentials on the line segment
$x\in[-1,1]$ correspond to the Chebyshev polynomials $\mT_k(x)$. 
From top to bottom: $\sigma_0(x)$, $\sigma_1(x)$, $\sigma_4(x)$.}
\label{fig:harmonic_basis}
\end{figure}

The corresponding harmonic functions $\phi^x_k$ are constructed via
\begin{equation}
\phi^x_k(\bx)=\int_{-H}^{H} 
\frac{-1\;}{2\pi}\ln\,|\bx-\xp\ex|\:\sigma_k(\xp;H)\,d\xp 
\label{eq:reconvolvex}
\end{equation}
and are illustrated in figure \ref{fig:harmonic_basis}.  
Specifying values along along the segment $[-1,1]$ in the $y$ direction, 
we obtain: 
\begin{equation}
\phi^y_l(\bx)=\int_{-1}^{1}
\frac{-1\;}{2\pi}\ln\,|\bx-y'\ey|\:\sigma_l(y';1)\,dy' 
\label{eq:reconvolvey}
\end{equation}
Note that $\sigma_k(x;H)=\sigma_k(x/H;1)/H$ for $k>0$ and
$\sigma_0(x;H)=\frac{\ln(1/2)}{H\ln(H/2)}\sigma_0(x/H;1)$.

The harmonic functions corresponding to specified values along the lower or
upper boundaries ($y=\pm 1$, $|x|< H$) are $\phi^x_k\left(\bx\mp
1\ey\right)$; those corresponding to the left or right boundaries ($x=\pm H$,
$|y|<1$) are $\phi^y_l\left(\bx\mp H\ex\right)$.  
We do not require the functions $\phi^x_k$, $\phi^y_l$ either inside or
outside $\Omega$, but only the values and normal derivatives on the boundary.
Although, for example, the values of $\phi^x_k\left(\bx + 1\ey\right)$ on the
lower boundary are merely the specified values $f_k(x)$, its values on the
other three boundaries must be calculated via \eqref{eq:reconvolvex}.  When
evaluating the normal derivatives, the kernel $G(\bx;x')$ is
differentiated before integration:
\begin{eqnarray}
\pd_y\phi^x_k(x,y)&=&\partial_y\int_{-H}^{H} 
\frac{-1\;}{2\pi}\ln\,\sqrt{(x-x')^2+y^2}\;\sigma_k(x';H)\,dx'\nonumber\\
&=&\int_{-H}^{H}\frac{-1\;}{2\pi}\frac{y}{(x-x')^2+y^2}\;\sigma_k(x';H)\,dx'
\label{eq:normalder}
\end{eqnarray}
and similarly for $\pd_x\phi^y_l$.

Any exterior harmonic function can be approximated 
by the truncated series:
\begin{eqnarray}
\phi^{K,L}(\bx)=
\sum_{k=0}^{K-1}\left[c^{x,-}_k\phi^x_k\left(\bx+H\ex\right) +
 c^{x,+}_k\phi^x_k\left(\bx-H\ex\right)\right]   \nonumber\\
+ \sum_{l=0}^{L-1}\left[c^{y,-}_l\phi^y_l\left(\bx+1\ey\right) 
+ c^{y,+}_l\phi^y_l\left(\bx-1\ey\right)\right]
\label{eq:truncser}\end{eqnarray}
The potential $\phi(\bx)$ of \eqref{eq:truncser} is defined by the $2(K+L)$
coefficients $\{c^{x,-}_k,c^{x,+}_k,c^{y,-}_l,c^{y,+}_l\}$.  A very important
property of the harmonic basis $\{\phi^x_k(\bx),\phi^y_l(\bx)\}$ is that it
represents a near-boundary field uniformly. This means that the truncated
series \eqref{eq:truncser} converges uniformly for any smooth boundary data
and for all locations $\bx$ near the boundary. This is a direct consequence of
the excellent convergence properties of Chebyshev approximation applied to
\eqref{eq:Carleman_equation}; a proof can be found in \cite{Levesley}.  This
property does not necessarily hold for other harmonic bases, in particular
spherical harmonics, for which near-boundary convergence cannot be achieved,
leading to a strong Gibbs effect.

In the taxonomy of boundary integral methods, equations
\eqref{eq:chargedist} and \eqref{eq:phiconstruct} constitute an indirect
method, in that the intermediate surface charge density $\sigma$ is
constructed; this is done by solving the
Fredholm integral equation of the first kind \eqref{eq:Carleman_equation}.
The surface charge density is a single-layer rather than a double-layer
(dipole) potential; equivalently $G$, rather than $\pd G/\pd n$, is used in
the representation.  Because the method determines only the boundary values
and normal derivatives of the exterior solution, it is not the preferred
approach when the exterior potential is itself required at each time step:
although the exterior solution can be sampled at any location, this is
computationally expensive, as is often the case for boundary integral methods.

\section{Electrostatic example}

We apply our method to a simple problem from electrostatics, 
the distribution of electric charges $-\rho_m$ 
confined in a rectangular domain but localized around the origin:
\begin{equation}
\rho_m(r,\theta)=\left\{\begin{array}{ll}
r^me^{-r^2/\delta^2}\cos(m\theta) &\text{for }|x|\leq 1\text{ and }|y|\leq 1 \\
0 &\text{for } |x|>1 \text{ or } |y|>1 \end{array}\right.
\label{eq:src_dist}
\end{equation}
The $r^m$ factor in (\ref{eq:src_dist}) ensures regularity of $\rho_m$ at
$r=0$.  The potential $\tilde{\Phi}_m$ due to {\it unbounded} sources (not
restricted to the interior domain) can be found analytically:
\begin{subequations}
\begin{eqnarray}
\tilde{\Phi}_{m=0}(r,\theta) &=& \frac{\delta^2}{4}\left[Ei\left(1,\frac{r^2}{\delta^2}\right)+2\log(r)\right]\label{eq:pot_sols_a} \\
\tilde{\Phi}_{m=1}(r,\theta) &=& \frac{\delta^4}{4r}\left[e^{-\frac{r^2}{\delta^2}}-1\right]\cos\theta \label{eq:pot_sols_b}\\
\tilde{\Phi}_{m=2}(r,\theta) &=& \frac{\delta^4}{4r^2}\left[\left(\delta^2+r^2\right)e^{-\frac{r^2}{\delta^2}}-\delta^2\right]\cos2\theta \label{eq:pot_sols_c}
\end{eqnarray}
\label{eq:pot_sols}
\end{subequations}\\
where $Ei(a,\lambda)$ is the error function 
$Ei(a,\lambda)\equiv \int_1^\infty{e^{-\lambda r}r^{-a}\,dr}$ and 
$\tilde{\Phi}_m$ are chosen to be finite at $r=0$.

We seek the corresponding electric potential. The parameter
$\delta$ is chosen to make $\rho_m$ very small near the boundaries. We
expect the solution to be almost unaffected by the presence of
boundaries. The source distribution $\rho_{m=0}(r,\theta)$ should 
therefore lead to a potential which is almost axisymmetric. 
Figure \ref{fig:axi_sol} shows the potential obtained
numerically for $\delta^2=0.15$ using the spectral resolution
$N=8$ in both directions. The domain boundary is represented by a bold
square. The contours are almost perfectly circular, as should be the 
case for $\delta$ small, showing that the presence of the boundaries
has minimal effect.

\begin{psfrags}
\psfrag{Convergence}[l][][1.0]{}
\psfrag{logError}[l][][0.9]{$\log_{10}E_{0}(N)$}
\begin{figure}
\begin{minipage}{0.45\textwidth}
\centering
\includegraphics[width=6.5cm]
{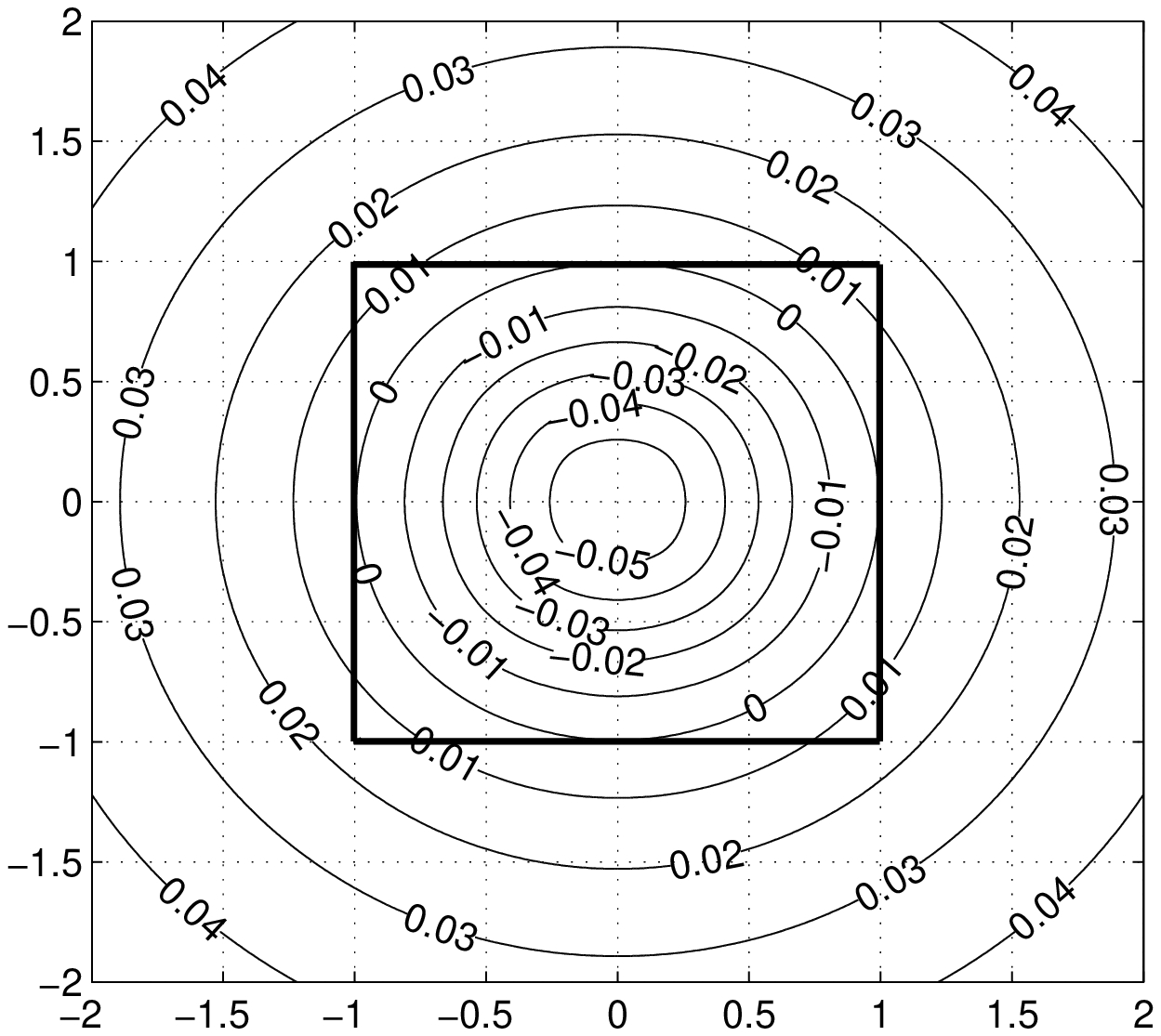}
\vspace{0.1cm}
\caption{Potential $\Phi_{m=0}^{N=8}$ for $\delta=0.15$. 
Maximal relative error is $E_{m=0}(N=8)\approx 0.03$.}
\label{fig:axi_sol}
\end{minipage}
\qquad
\begin{minipage}{0.45\textwidth}
\centering
\vspace{-0.2cm}
\includegraphics[width=6.9cm]
{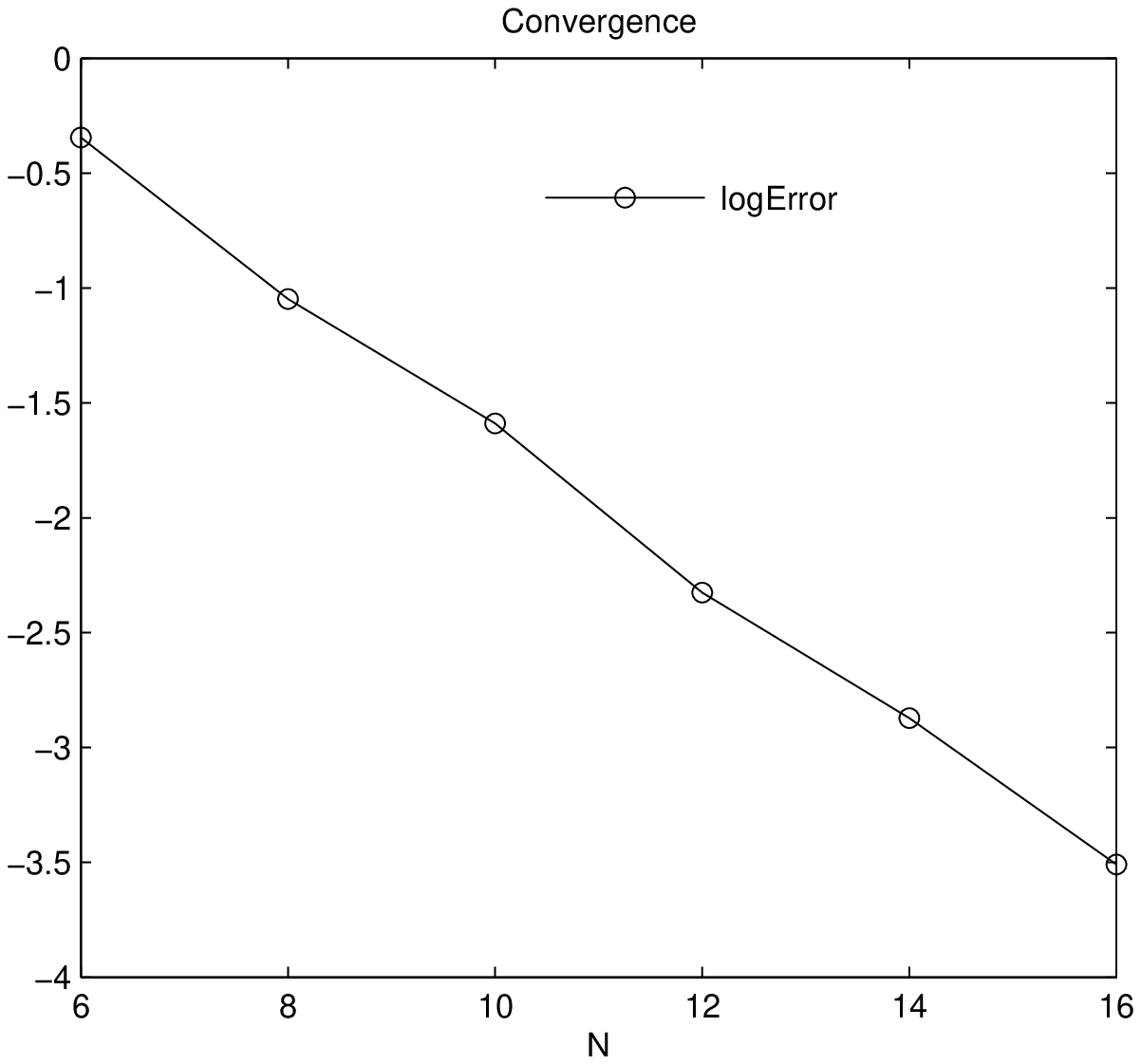}
\vspace{-0.1cm}
\caption{
Convergence test: Convergence test: $\log_{10}E_{m=0}(N)$ plotted for 
$N=[6,\ldots,16]$, $\delta^2=0.1$.}
\label{fig:conv}
\end{minipage}
\end{figure}
\end{psfrags}

To evaluate the error convergence of the method we computed the relative error
$E_m(N)$ defined as
\begin{equation}
E_m(N)=\sup_{r,\theta}\frac{|\tilde{\Phi}_m(r,\theta)-\Phi^N_m(r,\theta)|}{|\Phi_m(r,\theta)|}
\label{eq:error}
\end{equation}
where $\Phi^N_m(r,\theta)$ is the solution computed numerically with spectral
resolution $N$ in both spatial directions of the bounding square 
and $\tilde{\Phi}_m(r,\theta)$ is the analytic solution \eqref{eq:pot_sols} 
in the absence of the bounding square.
Figure \ref{fig:conv} proves the exponential convergence of the method.

\begin{figure}[htbp]
\begin{minipage}{0.45\textwidth}
\includegraphics[width=6.5cm]
{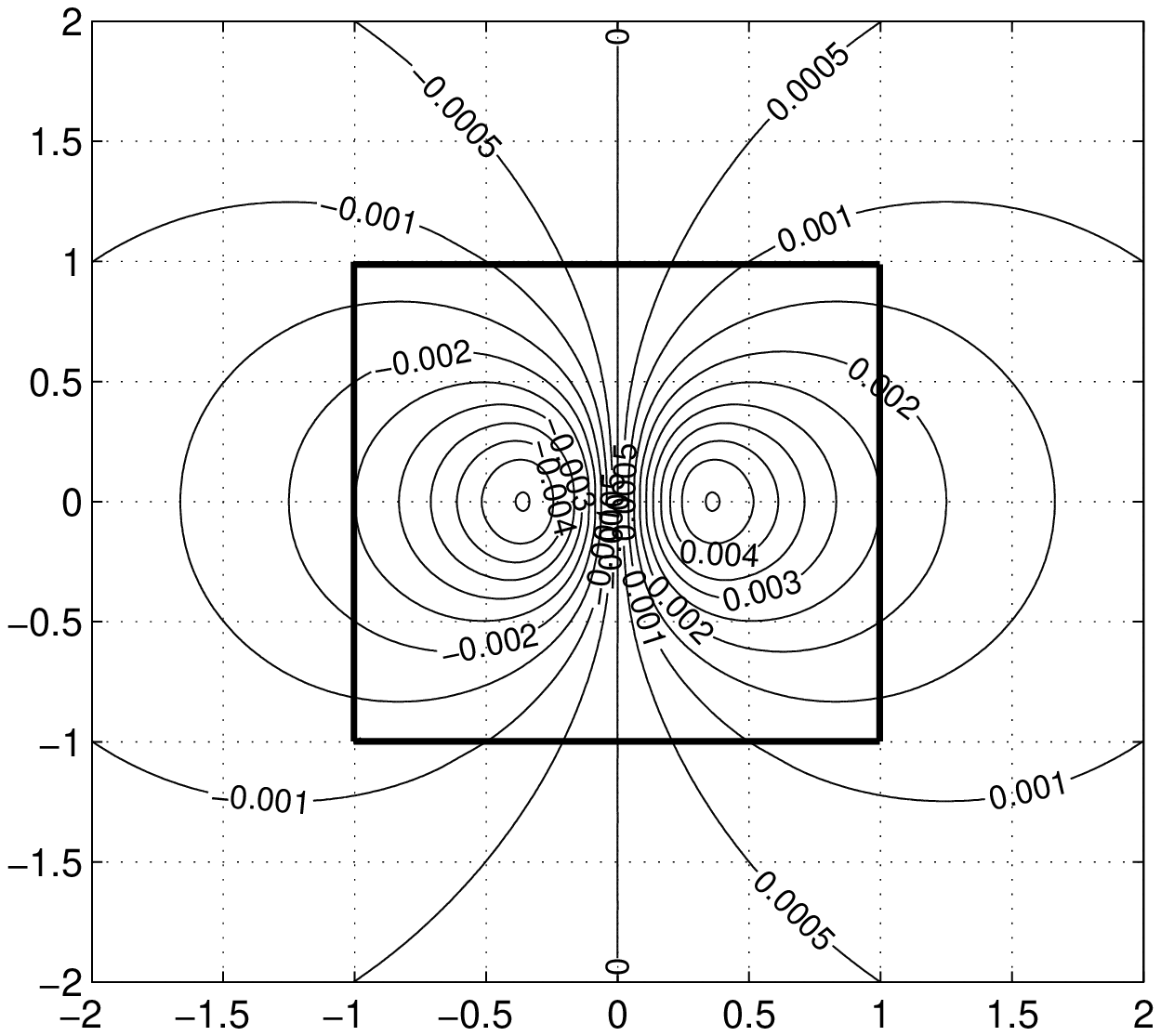}
\vspace*{.5cm}
\caption{Potential $\Phi_{m=1}^{N=16}$ for $\delta=0.1$.}
\label{fig:dipol}
\end{minipage}
\qquad
\begin{minipage}{0.45\textwidth}
\includegraphics[width=6.3cm]
{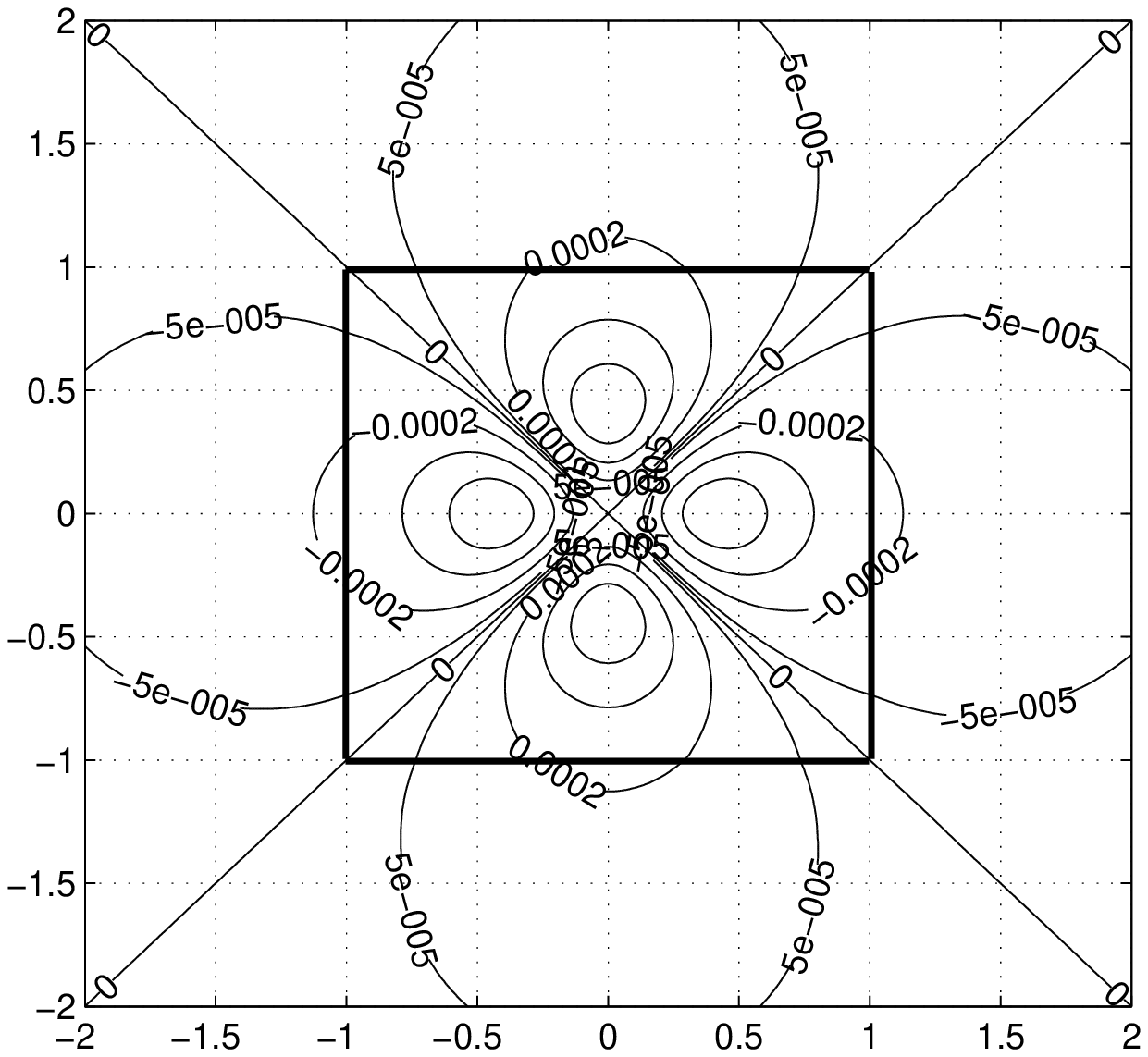}
\vspace*{.5cm}
\caption{Potential $\Phi_{m=2}^{N=16}$ for $\delta=0.1$.}
\label{fig:quad}
\end{minipage}

\vspace*{1cm}

\begin{minipage}{0.45\textwidth}
\includegraphics[width=6.5cm]
{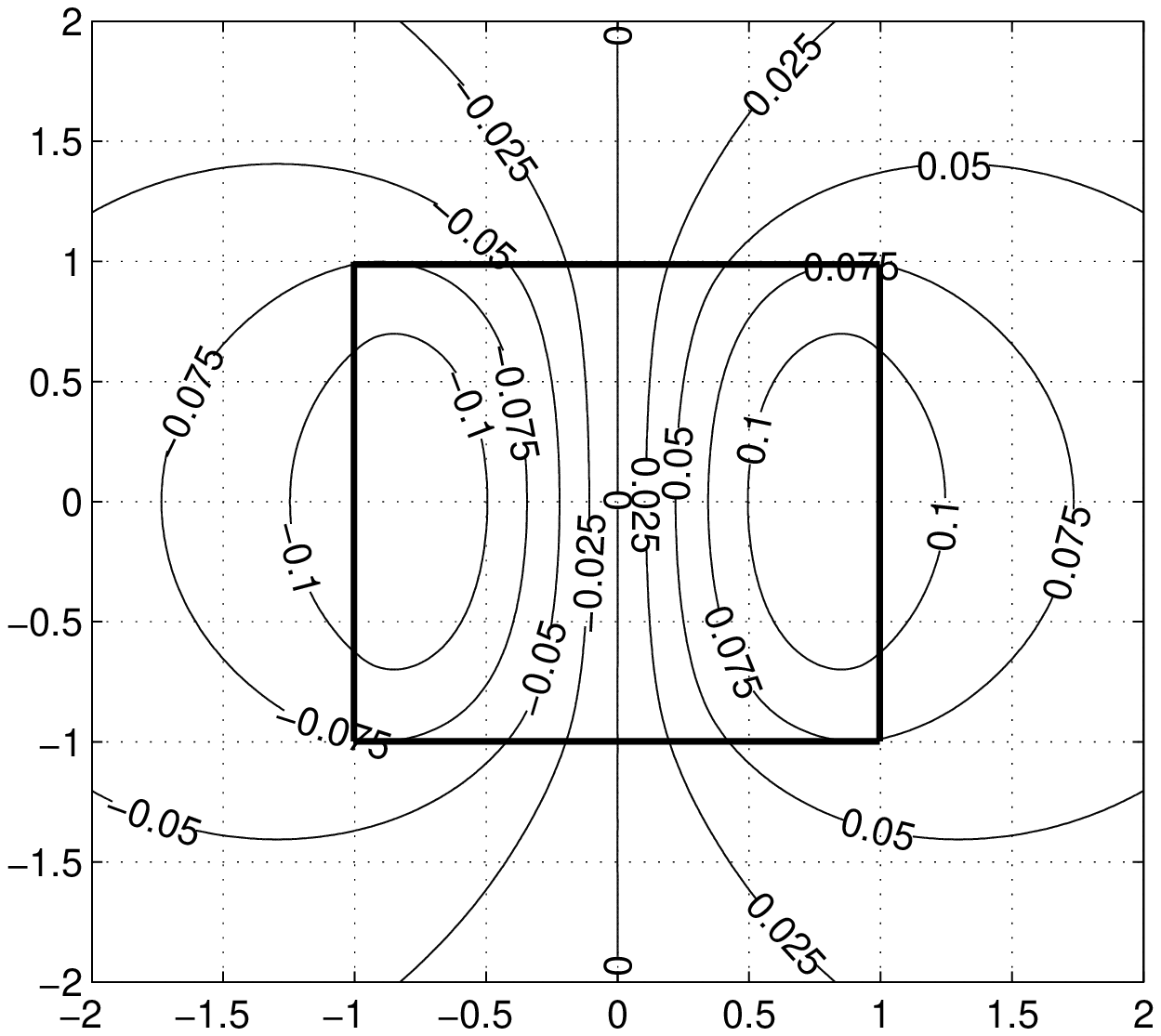}
\vspace*{.5cm}
\caption{Potential $\Phi_{m=1}^{N=16}$ for $\delta=2$.}
\label{fig:dipole}
\end{minipage}
\hspace*{.9cm}
\begin{minipage}{0.45\textwidth}
\includegraphics[width=6.5cm]
{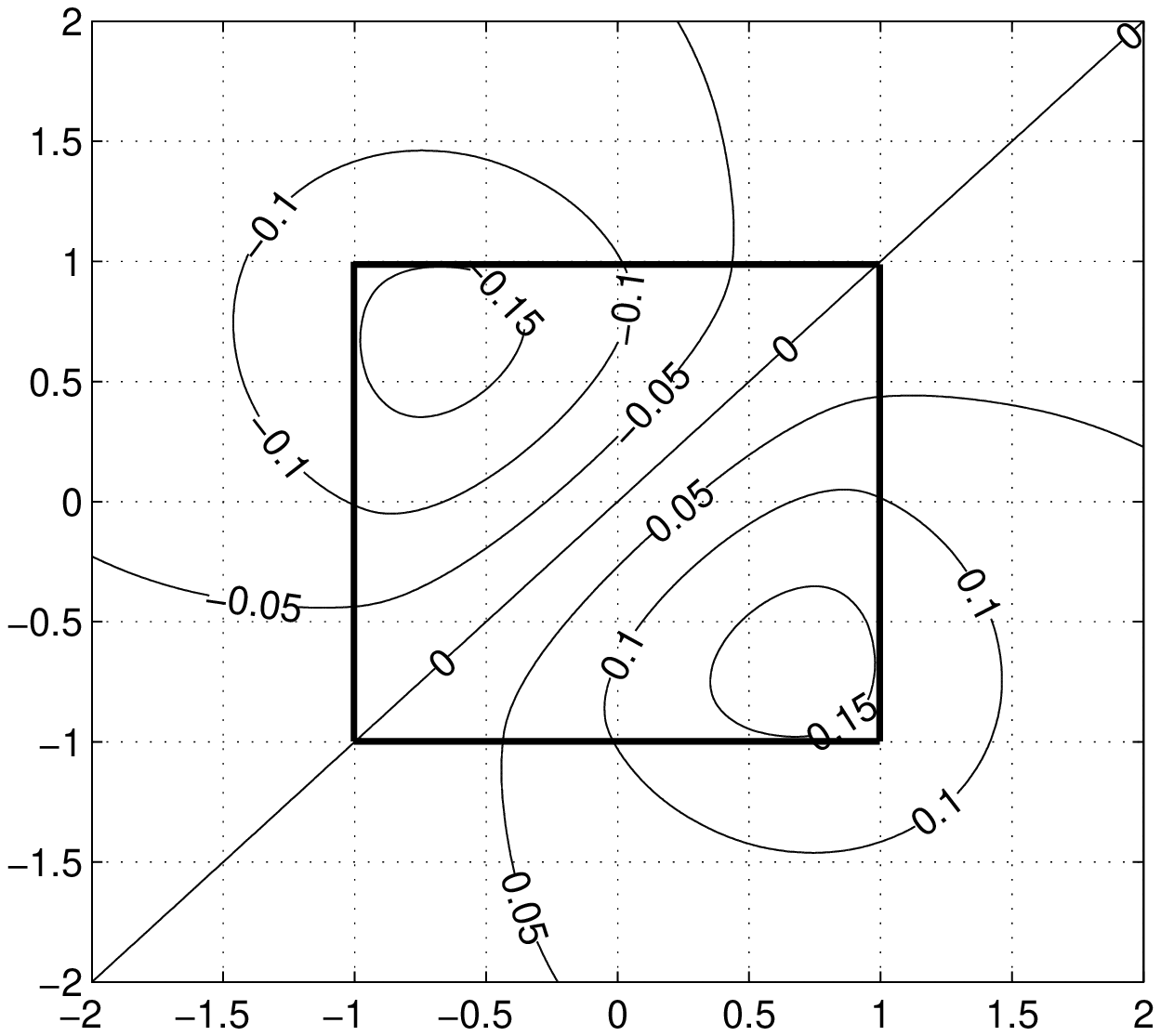}
\vspace*{.5cm}
\caption{$\Phi_{m=1}^{N=16}(r,\theta+\frac{\pi}{4})$ for $\delta=2$.}
\label{fig:turned_dipole}
\end{minipage}
\end{figure}

Figures \ref{fig:dipol}--\ref{fig:quad} show the electric potentials
$\Phi^{N=16}_{m=1}$ and $\Phi^{N=16}_{m=2}$ for $\delta=0.1$.  Figure
\ref{fig:dipole} shows $\Phi_{m=1}^{N=16}$ with $\delta=2$.  In figure
\ref{fig:turned_dipole}, the dipole source distribution has been rotated by
$45^\circ$ about the origin.  For this large value of $\delta$, charges are
located near the boundary.  In each case with $\delta \ll 1$, 
we observed exponential convergence toward solution \eqref{eq:pot_sols}. 
Convergence can only be confirmed up to a limited precision since the 
analytic solution \eqref{eq:pot_sols} does not
correspond exactly to the problem we are solving numerically, in which sources
are confined to the interior square.  The best agreement can be achieved for
small values of $\delta$.  If the numerical solution with highest spectral
resolution (here $N=64$) is instead taken as a reference, then the method
converges to this solution spectrally up to machine precision for any value 
of $\delta$.

\section{Implementation}

\subsection{Summary and computation cost}

We describe the implementation of the method for our illustrative example
\eqref{eq:ChebCheb} of the rectangle $[-H,H]\times [-1,1]$ with double
Chebyshev discretization $(K+1)\times (L+1)$.  

The total preprocessing step consists of: \\
$\bullet$ Evaluation of the values and the normal derivatives of
the exterior harmonic solutions on the boundary.\\
$\bullet$ Calculation of the interior harmonic solutions.\\
$\bullet$ Inversion or LU decomposition of the influence matrix, 
 
For each particular right-hand-side $\rho$, the operations consist of:\\
$\bullet$ Solving for a single particular solution.\\
$\bullet$ Acting with the inverse of the influence matrix.\\
$\bullet$ Using the corrected Dirichlet boundary conditions to
calculate the final solution.

The total number $J$ of
boundary points is $2(K+L)$.  The inversion or LU-decomposition of the
influence matrix $C$ in the preprocessing stage requires a time proportional
to $J^3$, while each solution of the linear system \eqref{eq:matchN3} 
determining the coefficients of the homogeneous solutions requires
a time proportional to $J^2$.
Each interior solution is calculated at a cost proportional to $KL^2$.

Symmetry can be used to reduce the cost of each step.
The symmetry of the rectangle divides all the independent 
harmonic solutions into four mutually orthogonal classes.
Decoupling the Laplacian operator according to parity
in $x$ and $y$ leads to four Poisson problems, each with resolution
$K/2 \times L/2$, thus reducing the time by a factor of two.
Decoupling by parity also reduces the influence matrix $C$ to four matrices,
the dimensions of which are one fourth of that of the original matrix.

Table \ref{tab:cost} gives the operation count of each step,
taking into account the reductions permitted by symmetry.
\begin{table}
\begin{tabular}{|l|ll|}
\cline{1-3}
Calculation &  Result & Cost  \\ \cline{1-3}
\multicolumn{3}{|l|}{Preprocessing} \\ \cline{1-3}
Exterior harmonic solutions & $\phi_j|_\bdy$, $\pdn\phi_j|_\bdy$ & 
$\begin{array}{l} K\times(\mbox{Int}_x + \mbox{Int}^\prime_x)\\
+L\times(\mbox{Int}_y + \mbox{Int}^\prime_y)
\end{array}$\\
Interior harmonic solutions & $\pdn\Phi^h_j$ & $(K+L)KL^2/4$ \\
Influence matrix inversion/decomposition & $C^{-1}$ & $(K+L)^3/2$ \\
\cline{1-3}
\multicolumn{3}{|l|}{For each right-hand-side} \\ \cline{1-3}
Particular solution & $\Phi^p$ & $KL^2/2$ \\
Action with influence matrix & $c_j$ & $(K+L)^2$ \\
Corrected solution & $\Phi$ & $KL^2/2$ \\\cline{1-3}
\end{tabular}
\label{tab:cost}
\caption{Operation count of each preprocessing and 
right-hand-side-dependent step for a rectangle
discretized with $(K+1)\times(L+1)$ Chebyshev polynomials and points.
$\mbox{Int}_x$ and $\mbox{Int}^\prime_x$ are the costs of performing
the singular integrals over $x$ in \eqref{eq:reconvolvex} 
and \eqref{eq:normalder}, and
$\mbox{Int}_y$ and $\mbox{Int}^\prime_y$ are those of
the analogous integrals over $y$.}
\end{table}\subsection{Singular integrals}

The integrations in \eqref{eq:reconvolvex}-\eqref{eq:normalder} are
performed numerically. Special attention must be paid in doing so since both
the kernel $G(\bx;x')$ and the density $\sigma(x')$ have integrable
singularities within the domain of integration. The singular points are
$x'\ex=\bx$ and $x'=\pm H$ for $\phi^x_k$ and $y'\ey=\bx$ and $y'=\pm 1$ for
$\phi^y_l$.  Dedicated adaptive quadratures (see \cite{Press86}) can be used to
compute these integrals accurately.

It is also possible to evaluate the singular part of the integral
analytically, reducing the numerical problem to the evaluation of integrals
with non-singular integrands.  The remaining integrand is piecewise $C^\infty$
and can be integrated with spectral precision over each of the regular
subdomains.  Singularity subtraction greatly decreases the variation in grid
density needed to sample the integrand homogeneously, thereby significantly
accelerating the numerical integrations in 
\eqref{eq:reconvolvex}-\eqref{eq:normalder}.
Specifically, an adaptive method requires a smaller
number of iterations, or, alternatively, a non-adaptive method requires a
coarser resolution.  However, the convergence of our approximation with $K,L$
is exponential (spectral), regardless of whether the singular part of the
integral is subtracted or included in the numerical evaluation.

\subsection{Conditioning of matrices}

\begin{psfrags}
\psfrag{fit}[r][][0.9]{$3.58N^2-18.16N+59$}
\psfrag{C(N)}[r][][0.9]{$\mathcal{C}(N)$}
\begin{figure}
\centering
\includegraphics[width=8cm]{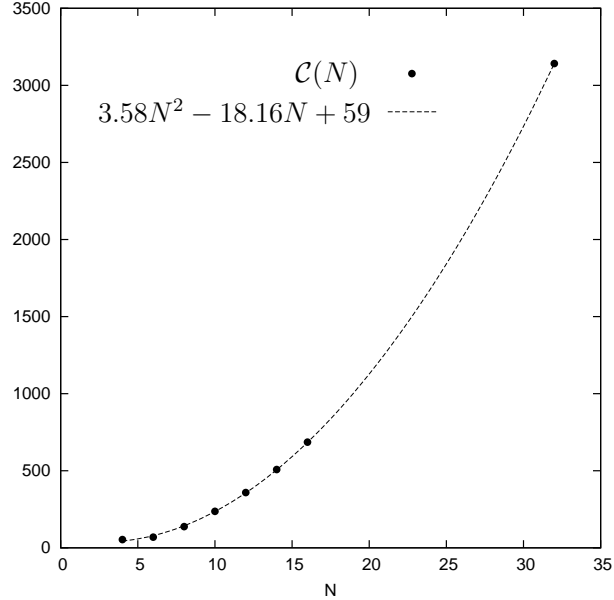}
\label{fig:cond_fit}
\caption{Quadratic fit of the condition number $\mathcal{C}(N)$ of the 
influence matrix defined in \eqref{eq:def_influ_matrix}. }
\end{figure}
\end{psfrags}

The influence matrix \eqref{eq:def_influ_matrix} is not immediately
invertible.  Because of the redundancy of information at the corners, this
matrix has exactly four zero singular values or eigenvalues. The corresponding
linear system can be solved after arbitrarily correcting singular values or
eigenvalues of the influence matrix; see
\cite{Tuckerman89,Boronski05} for more details.

The condition number $\mathcal{C}$ of the corrected matrix depends on the
spatial resolution $N$ and the maximal order of derivatives
used to express the boundary conditions.  In our case of a Neumann boundary
condition and a resolution $N$ in each direction,
the condition number scales as $\mathcal{C}=O\left(N^2\right)$.
Fitting the condition numbers computed for $N\in[2,32]$ with a parabola
(see fig. \ref{fig:cond_fit}) yields a formula for predicting the condition
number for an arbitrary resolution:
\begin{equation}
\mathcal{C}(N)=3.58N^2-18.16N+59
\label{eq:cond_number}
\end{equation}
It can then be deduced from \eqref{eq:cond_number} that a reasonably
conditioned matrix with $\mathcal{C}<10^7$ is obtained for
a spatial resolution as high as $O(1000)$.
The method can therefore be applied to problems where
small-scale field features require use of high spatial resolution.

\section{Generalizations}

We now discuss the applicability of this method to other geometries,
problems, and spatial discretizations.
The decomposition into interior particular and homogeneous functions
and exterior homogeneous functions described in section 
\ref{sec:Influence matrix formulation} is, of course, 
completely general and not related to any
particular spatial discretization.
The method described in section \ref{sec:Interior domain}
for constructing the exterior harmonic functions
relies on the Chebyshev-Chebyshev discretization of the rectangle
which is widely used since the Chebyshev 
polynomials are optimal approximants of smooth functions. 
This property, as well as the straightfoward correspondence
between interior and exterior solutions, 
make the Chebyshev discretization especially suitable
for the construction of the exterior harmonic solutions as well.
However, the method is easily generalizable to other basis functions 
$f_k(x)$ for the potential values, which can be
substituted into \eqref{eq:Carleman_sol} in order to calculate the 
corresponding charge densities $\sigma_k(x)$, if dictated
by the geometry or numerical method used for the interior problem.

In fact, since our real interest is in generating the complete set
of $\sigma_k(x)$ necessary to generate the complete set of $\phi_k(x)$, 
rather than in calculating the specific $\sigma_k(x)$
corresponding to each particular $f_k(x)$,
the only information really required in \eqref{eq:Carleman_sol_Cheb} 
is the singularity 
$1/\pi\sqrt{1-x^2}$. One may then allow the set of $\sigma_k$'s to be the 
products of this singularity with the members of any appropriate basis set 
of analytic functions on the boundary in question.

In three dimensions, the fundamental Green's function is
\begin{equation}
G(\bx,\bxp) = \frac{1\;}{4\pi}\frac{1}{|\bx-\bxp|}
\end{equation}
In an axisymmetric geometry with a Fourier representation of the
azimuthal direction, all of the problems to be solved decouple according to
Fourier mode. The operation count would then scale linearly with the number of
azimuthal points or Fourier modes.  The elliptic problems in
\eqref{eq:construct} would remain two-dimensional, and the integral equations
equivalent to \eqref{eq:Carleman_equation} would remain one-dimensional.  

This method can also be applied to other elliptic problems
or to parabolic problems. 
As stated in section \ref{sec:motivation}, our motivation
for developing this method is to apply it to the magnetohydrodynamic
equations \eqref{eq:Bint}-\eqref{eq:Bcont}, in which
\eqref{eq:Bint} is a parabolic equation.
A general parabolic problem may be written as:
\begin{equation}
\pd_t \Phi = \lap \Phi + \mathcal{F}(\Phi,\rho)
\label{eq:parabolic}\end{equation}
where $\mathcal{F}$ may include nonlinear and/or time-dependent source terms.
First-order implicit temporal discretization of \eqref{eq:parabolic}
results in the inhomogeneous Helmholtz equation:
\begin{equation}
(I-\delta t \lap) \Phi(t+\delta t) = \mathcal{F}
\label{eq:implicit}\end{equation}
where $\mathcal{F}$ may depend on previous values of $\Phi$.
This Helmholtz operator $(I-\delta t \lap)$ can replace the Laplacian in 
\eqref{eq:two_domains} and \eqref{eq:construct}.
It is known that replacement of the Helmholtz equation by a boundary 
integral equation can lead to singularities for certain values of the 
wavenumber (here $i/\sqrt{\delta t}$); a large body of work, 
e.g.~\cite{Amini95,Givoli98,Gerdes96,Hwang99,TCLin04}, addresses this problem.
However, in the magnetohydrodynamic case of a conducting fluid surrounded
by an exterior vacuum, no such difficulties would be introduced,
since the exterior problem remains governed by Laplace's equation.
More complicated vectorial operators may appear,
as occur in the Navier-Stokes or magnetohydrodynamic equations.

\section{Comparison with other approaches}

We mention here some other techniques that have been used
to solve exterior problems or to match interior and exterior domains.
Spectral methods can be combined with
various transformations and mappings.
The inner region can be surrounded by a sphere, and the
outer domain decomposed into the region inside and outside the sphere.
The exterior domain can be mapped into an interior domain via
a $1/r$ mapping \cite{Jepps67,Grandclement01,Lai06};
spectral methods can then be used to treat either or both domains.
The region exterior to one or more spheres has been mapped to the
interior of a rectangle \cite{Ansorg04} or a pentangular
\cite{Ansorg05} region rotated about an axis, and Chebyshev-Fourier
expansions used to solve elliptic equations arising in the study
of black holes in general relativity.
A smooth boundary can be parameterized by angle, and the boundary
values represented as a series of trigonometric functions
or spherical harmonics
\cite{Atkinson82,Gerdes96,Hwang99,TCLin04,Grandclement01,Lai06,Meddahi02,Chen94,Ganesh98}.
Our method differs from these in that
a Chebyshev approximation is used to represent
the boundary values on each segment of a non-smooth boundary,
and an analytic formula is used to calculate the surface density
which exactly yields this Chebyshev approximation.

Conformal mapping is another technique which can be used to calculate 
interior or exterior harmonic functions.
The Riemann mapping theorem guarantees the existence of a conformal
transformation from the interior or exterior of a simply connected domain into
the interior or exterior of a unit disk; its proof is, however,
non-constructive, and does not explicitly derive the transformation.  For
some geometries, including the exterior of a rectangle, an analytical formula
can be derived \cite{Ivanov}.  For polygon-bounded regions
with piecewise-constant boundary conditions, the
Schwarz-Christoffel \cite{TrefethenSC} mapping has proved to be a very robust
tool, applied to problems arising in magneto- and electro- statics, potential
flows, inverse problems and many other fields.  

Our influence matrix approach relies on calculating
harmonic functions with arbitrary boundary data,
for which conformal mapping is much more problematic.
More general conformal mappings are often computed by 
solving Symm's or Carleman's equation \eqref{eq:Carleman_equation} 
numerically on
the domain boundary, making this approach similar in terms of 
numerical cost and precision to boundary integral equation methods.
It is interesting to note that, for domains including
corners, the Chebyshev approximation is especially well suited, guaranteeing 
superconvergence of the mapping function \cite{Levesley}.

\section{Conclusion}

As a test case for the magnetohydrodynamic equations, we have developed a
method for solving the two-dimensional Poisson equation in a bounded domain,
where the solution satisfies matching conditions with a harmonic potential outside the
domain. The method solves only the interior problem and determines the
boundary conditions ensuring smooth matching with the exterior solution. The
essential element of this approach is construction of a basis of harmonic
functions which represent the near-boundary exterior solutions uniformly.
This basis is used to construct the influence matrix which serves to
impose the coupled boundary conditions between the interior and exterior
solutions. The method is numerically reasonably well conditioned and can be
used for quite high spatial resolutions. For a spectral solver, this method
guarantees exponential convergence.

Instead of corresponding to point sources on the boundary, each exterior
harmonic solution corresponds to a spectral basis function.  The most costly
process -- the construction of a basis of exterior harmonic functions --
depends only on geometry and spatial resolution.  Once the basis is computed
it can be stored and used for any computation using the same resolution and
domain shape.  When used as a preprocessing step for time-dependent
simulations, the cost of constructing the exterior harmonic basis is
negligible compared to that of thousands of time steps.  

Since cylindrical coordinates have one periodic direction, it should
be possible to apply this method separately to each of the Fourier modes,
treated individually as two-dimensional problems.
The extension of this method to the MHD equations in a finite 
cylindrical geometry is currently under investigation.


\begin{thebibliography}{99}

\bibitem{Steglitz02} 
{\sc R.~Steglitz \& U.~M\"uller}, 
Experimental demonstration of a homogenous two-scale dynamo, 
Magnetohydrodynamics~{\bf 38}, 27 (2002).

\bibitem{Riga01}
{\sc A.~Gailitis, O.~Lielaisis, E.~Platacis, S.~Dementev, A.~Cifersons,
  G.~Gerbeth, T.~Gundrum, F.~Stefani, M.~Christen \& G.~Will,}
Magnetic field saturation in the Riga dynamo experiment,
Phys. Rev. Lett.~{\bf 86}, 3024 (2001).

\bibitem{VKS02}
{\sc M.~Bourgoin, L.~Mari\'e, F.~Petrelis, C.~Gasquet, A.~Guigon, J.B.~Luciani,
  M.~Mulin, F.~Namer, J.~Burgete, A.~Chiffaudel, F.~Daviaud, S.~Fauve,
  P.~Odier \& J.F.~Pinton}
Magnetohydrodynamics measurements in the von K\'arm\'an sodium experiment,
Phys. Fluids~{\bf 14}, 3046 (2002).

\bibitem{Forest02}
{\sc C.B.~Forest, R.A.~Bayliss, R.D.~Kendrick, M.D.~Nornberg,
  R.~O'Connell \& E.J.~Spence,}
Hydrodynamic and numerical modeling of a spherical homogeneous dynamo 
experiment,
Magnetohydrodynamics~{\bf 38}, 107 (2002).

\bibitem{Shew02}
{\sc  W.L.~Shew, D.R.~Sisan, \& D.P.~Lathrop,}
Mechanically forced and thermally driven flows in liquid sodium,
Magnetohydrodynamics~{\bf 38}, 121 (2002).

\bibitem{Jepps67}
{\sc S.A.~Jepps}, 
Numerical models of hydromagnetic dynamos,
J.~Fluid Mech.~{\bf 67}, 625 (1967).

\bibitem{Dudley89}
{\sc M.~Dudley \& R.~James}, Time-dependent kinematic dynamos with
stationary flows, Proc. Roy. Soc. London A {\bf 425}, 407--429
(1989).

\bibitem{Glatzmaier95}
{\sc G.~A.~Glatzmaier and P.~H.~Roberts}
A three-dimensional self-consistent computer simulation of a geomagnetic field reversal,
Nature~{\bf 337}, 203 (1995).

\bibitem{Tilgner97}
{\sc A.~Tilgner},
A kinematic dynamo with a small scale velocity field,
Phys.~Rev.~{\bf A 226}, 75--79 (1997).

\bibitem{Hollerbach00}
{\sc R.~Hollerbach}, A spectral solution of the magneto-convection
equations in spherical geometry,
Int.~J.~for Num.~Meth.~in Fluids {\bf 32}, 773--797 (2000).

\bibitem{Willis02}
{\sc A.~P.~Willis \& C.F.~Barenghi,}
A Taylor-Couette Dynamo,
Astronomy \& Astrophysics~{\bf 393}, 339--343 (2002).

\bibitem{Matsui04}
{\sc H.~Matsui \& H.~Okuda},
Development of a simulation code for MHD dynamo processes using the GeoFEM
platform, Intl.~J.~of Comp.~Fluid Dyn.~{\bf 18}, 323--332 (2004).

\bibitem{Iskakov04}
{\sc A.~B.~Iskakov, S.~Descombes \& E.~Dormy},
An integro-differential formulation for magnetic induction in
bounded domains: boundary element-finite volume method,
J.~Comput.~Phys.~{\bf 7}, 540--554 (2004).

\bibitem{Iskakov05}
{\sc A.~B.~Iskakov \& E.~Dormy},
On magnetic boundary conditions for non-spectral dynamo simulations,
Geophys.~Astrophys.~Fluid Dyn.~{\bf 99}, 481--492 (2005).
  
\bibitem{Xu04}
{\sc M.~Xu, F.~Stefani \& G.~Gerbeth},
The integral equation method for a steady kinematic dynamo problem,
J.~Comput.~Phys.~{\bf 196}, 102--125 (2004).

\bibitem{Guermond06}
{\sc J.L.~Guermond, R.~Laguerre, J.~L\'eorat \& C.~Nore},
An interior penalty Galerkin method for the MHD equations in
heterogeneous domains, J. Comput. Phys.~{\bf 221}, 349--369 (2007). 

\bibitem{Roberts67}
{\sc P.~Roberts}, An Introduction to Magnetohydrodynamics,
Longmans Green, London, 1967.

\bibitem{Hsiao73}
{\sc G.~Hsiao \& R.C.~Maccamy}
Solution of Boundary Value Problems by Integral Equations of the First Kind
SIAM Review {\bf 15}, 687--705 (1973).

\bibitem{Jaswon77}
{\sc M.A.~Jaswon \& G.T.~Symm},
Integral equation methods in potential theory and elastostatics,
Academic Press, 1977.

\bibitem{Brebbia78}
{\sc C.~Brebbia \& S.~Walker},
The boundary element method for engineers,
Pentech Press, 1978.

\bibitem{Banerjee81}
{\sc P.K.~Banerjee \& R.~Butterfield},
Boundary element methods in engineering science,
McGraw-Hill, 1981.

\bibitem{Atkinson82}
{\sc K.~E.~Atkinson},
The numerical solution of Laplace's equation in three dimensions,
SIAM J.~Num.~Anal.~{\bf 19}, 263--274 (82).

\bibitem{Rokhlin83}
{\sc V.~Rokhlin},
Rapid solution of integral equations of classical potential theory,
J.~Comput.~Phys. {\bf 60}, 187--207 (1983).

\bibitem{Brebbia83}
{\sc C.A.~Brebbia, T.~Futagami \& M.~Tanaka, eds.},
Boundary Elements: Proceedings of the 5th International Conference,
Springer, 1983.

\bibitem{Pozrikidis92}
{\sc C.~Pozrikidis},
Boundary integral and singularity methods for linearized viscous flow,
Cambridge Univ. Press, 1992.

\bibitem{KirkupWeb}
{\sc S.~Kirkup},
{\tt http://www.boundary-element-method.com}

\bibitem{Orszag71}
{\sc S.A.~Orszag}, Numerical simulation of incompressible flows within
simple boundaries: I. Galerkin (spectral) representations,
Stud.~Appl.~Math.~{\bf 50}, 293--327 (1971).

\bibitem{Canuto88}
{\sc C.~Canuto, M.Y.~Hussaini, A.~Quateroni \& T.A.~Zang},
Spectral Methods in Fluid Dynamics, Springer, 1988.

\bibitem{Jorgens70}
{\sc K.~J\"orgens},
Integral operators, Pitman, Boston, 1982.
Translation from German of Lineare Integraloperatoren,
Teubner, Stuttgart, 1970.

\bibitem{Polyanin98}  
{\sc A.~D.~Polyanin \& A.~V.~Manzhirov},
Handbook of Integral Equations, CRC Press, Boca Raton, 1998.
  
\bibitem{Mason00}
{\sc J.C.~Mason \& D.C.~Handscomb},
Chebyshev Polynomials, Chapman \& Hall/CRC Press, 2000.

\bibitem{Trefethen00}
{\sc L.~N.~Trefethen}, Spectral Methods in Matlab (SIAM, Philadelphia, 2000).

\bibitem{Levesley} 
{\sc J.~Levesley, D.~M.~Hough \& S.N.~Chandler-Wilde}, 
A Chebyshev collocation method for solving Symm's integral equation 
for conformal mapping: a partial error analysis, 
IMA J.~Num.~Anal. {\bf 14}, 57--79 (1993).

\bibitem{Press86}
{\sc W.H.~Press, B.P.~Flannery, S.A.~Teukolsky \& W.T.~Vetterling},
Numerical Recipes: The Art of Scientific Computing, Cambridge
Univ.~Press, 1986.

\bibitem{Tuckerman89}
{\sc L.S.~Tuckerman}, 
Divergence-free velocity fields in nonperiodic geometries, 
J.~Comput.~Phys.~{\bf 80}, 403--441 (1989).

\bibitem{Boronski05}
{\sc P.~Boronski}, 
A Method Based on Poloidal-Toroidal Potentials Applied to the von K\'arm\'an 
Flow in a Finite Cylinder Geometry,
Ph.D.~Thesis, Ecole Polytechnique, 2005.

\bibitem{Amini95}
{\sc S.~Amini \& S.M.~Kirkup},
Solution of Helmholtz equation in the exterior domain by
elementary boundary integral methods,
J.~Comput.~Phys.~{\bf 118}, 208--221 (1995).

\bibitem{Givoli98}
{\sc D.~Givoli \& I.~Harari}, eds.,
Special Issue on Exterior Problems of Wave Propagation,
Comput.~Methods Appl.~Mech.~Engrg.~{\bf 164}, 1-266 (1998).

\bibitem{Gerdes96}
{\sc K.~Gerdes \& L.~Demkowicz},
Solution of 3D-Laplace and Helmholtz equations in
exterior domains using hp-infinite elements,
Comput.~Methods Appl.~Mech.~Engrg.~{\bf 137}, 239--272 (1996).

\bibitem{Hwang99}
{\sc W.S.~Hwang},
A boundary spectral method for solving exterior
acoustical problems with hypersingular integrals,
Int.~J.~Numer.~Meth.~Engng.~{\bf 44}, 1775--1783 (1999).

\bibitem{TCLin04}
{\sc T.-C.~Lin \& Y.~Warnapala-Yehiya},
The numerical solution of exterior Neumann problem
for Helmholtz's equation via modified Green's functions approach,
Computers and Mathematics with Applications {\bf 47}, 593--609 (2004).

\bibitem{Grandclement01}
{\sc P.~Grandclement, S.~Bonazzola, E.~Gourgoulhon \& J.-A.~Marck},
A multidomain spectral method for scalar and
vectorial Poisson equations with noncompact sources,
J.~Comput.~Phys.~{\bf 170}, 231--260 (2001).

\bibitem{Lai06}
{\sc M.-C.~Lai, Z.~Li \& X.~Lin},
Fast solvers for 3D Poisson equations involving interfaces
in a finite or the infinite domain,
J.~Comput.~Appl.~Math.~{\bf 181}, 106--125 (2006).

\bibitem{Ansorg04}
{\sc M.~Ansorg, B.~Br\"ugmann \& W.~Tichy},
A single-domain spectral method for black hole puncture data,
Phys.~Rev.~D {\bf 70}, 064011 (2004).

\bibitem{Ansorg05}
{\sc M.~ Ansorg},
Double-domain spectral method for black hole excision data,
Phys.~Rev.~D {\bf 72}, 024018 (2005).

\bibitem{Meddahi02}
{\sc S.~Meddahi \& A.~M\'arquez},
A combination of spectral and finite elements
for an exterior problem on the plane,
Appl.~Numer.~Math.~{\bf 43}, 275--295 (2002).

\bibitem{Chen94}
{\sc Y.~Chen},
Galerkin methods for solving single layer integral equations
in three dimensions, Ph.D.~Thesis, Univ.~of Iowa, 1994.

\bibitem{Ganesh98}
{\sc M.~Ganesh, I.G.~Graham \& J.~Sivaloganathan},
A new spectral boundary integral collocation method
for three-dimensional potential problems,
SIAM J.~Num.~Anal.~{\bf 35}, 778--805 (1998).

\bibitem{Ivanov}
{\sc V.I.~Ivanov \& M.K.~Trubetskov}, 
Handbook of Conformal Mapping with Computer-Aided Visualization, CRC Press, Boca Raton, 1995.

\bibitem{TrefethenSC}
{\sc T.A.~Driscoll \& L.N.~Trefethen}, Schwarz-Christoffel Mapping, 
Cambridge Univ.~Press, 2002.

\end{thebibliography}
\end{document}